\documentclass[letterpaper, 10 pt, conference]{ieeeconf}  % Comment this line out
\IEEEoverridecommandlockouts                              % This command is only
\usepackage[english]{babel}

\newtheorem{theorem}{Theorem}
\newtheorem{lemma}{Lemma}
\newtheorem{remark}{Remark}
\newtheorem{assumption}{Assumption}
\newtheorem{corollary}{Corollary}

\usepackage{graphics} % for pdf, bitmapped graphics files
\usepackage{epsfig} % for postscript graphics files
\usepackage{times} % assumes new font selection scheme installed
\usepackage{amsmath}

\usepackage{amssymb}
\usepackage{physics}
\usepackage{mathtools}

\usepackage{algorithm}
\usepackage{algpseudocode} 

\usepackage{multicol}

\usepackage{subcaption}
\usepackage{caption}

\usepackage{xcolor}

\def\mz{\mu_z}
\def\mt{\mu_\theta}
\def\va{v_\alpha}
\def\vas{v_\alpha^\star}

\newcommand{\ceil}[1]{\left\lceil {#1} \right\rceil}

\newcommand\scalemath[2]{\scalebox{#1}{\mbox{\ensuremath{\displaystyle #2}}}}

\newcommand*{\QEDB}{\hfill\ensuremath{\blacksquare}}%

\newcommand*{\vsepfbox}[1]{%
  \begingroup
    \sbox0{\fbox{#1}}%
    \setlength{\fboxrule}{0pt}%
    \mbox{\kern-\fboxsep\fbox{\unhbox0}\kern-\fboxsep}%
  \endgroup
}
\newcommand{\hide}[1]{}
\title{\LARGE \bf
Distributed Primal-dual Optimization for Heterogeneous Multi-agent Systems
}

\author{Yichuan Li$^{1}$, Petros Voulgaris$^{2}$, and Nikolaos M. Freris$^{3}$ % <-this % stops a space
\thanks{$^{1}$Coordinated Science Laboratory, University of Illinois at Urbana-Champaign, IL 61820, USA. {\tt\small yli129@illinois.edu}. }  %
\thanks{$^{2}$Department of Mechanical Engineering, University of Nevada, Reno, NV 89557, USA. {\tt\small pvoulgaris@unr.edu}.}
\thanks{$^{3}$ School of Computer Science, University of Science and Technology of China, Hefei, Anhui, 230027, China. {\tt\small nfr@ustc.edu.cn.}}
\thanks{*This work was supported by the Ministry of Science and Technology of China under grant 2019YFB2102200. Correspondence to N. Freris.}
}

\begin{document}
\maketitle

%%%%%%%%%%%%%%%%%%%%%%%%%%%%%%%%%%%%%%%%%%%%%%%%%%%%%%%%%%%%%%%%%%%%%%%%%%%%%%%%
\begin{abstract}
Heterogeneous networks comprise agents with varying capabilities in terms of computation, storage, and communication. In such settings, it is crucial to factor in the operating characteristics in 
allowing agents to choose appropriate updating schemes, so as to better distribute computational tasks and utilize the network more efficiently. 
We consider the multi-agent optimization problem of cooperatively minimizing the sum of local strongly convex objectives. We propose an asynchronous distributed primal-dual protocol, which allows for the primal update steps to be agent-dependent (an agent can opt between first-order or Newton updates). Our analysis introduces a unifying framework for such hybrid optimization scheme and establishes global linear convergence  in expectation, under strongly convex objectives and general agent activation schemes. Numerical experiments on real
life datasets attest to the merits of the proposed algorithm.
\end{abstract}

\section{INTRODUCTION}
By harnessing the computational power of multiple machines, distributed optimization aims to integrate a network of agents to solve a minimization problem. The framework has found success across a wide range of fields, including control, smart-grid networks, and machine learning \cite{angelia2018}--\nocite{huang2020}\nocite{babi}\cite{khan}. The growing number of mobile devices with on-device computation and communication capabilities further instigates the possibility of ever-larger scale networks with decentralized computing, such as Cyber-Physical Systems (CPS) \cite{CPS2012}. In this setting, heterogeneity is ubiquitous and relates to several aspects of the system such as: hardware conditions, network connectivity, data volumes, as well as communication bandwidths. These challenges require revisiting algorithmic designs to account for the operating characteristics to effectively distribute computational tasks. 

We consider the following collaborative network problem: 
\begin{align}
    \underset{\Tilde{x}\in\mathbb{R}^d}{\mathrm{minimize}}\,\,\left\{ \sum_{i=1}^m f_i(\Tilde{x})+g(\Tilde{x})\right\}, \label{prob1}
\end{align}
where $f_i(\cdot):\mathbb{R}^d\to \mathbb{R}$ are local strongly convex objectives, and $g(\cdot)$ is a convex but possibly nonsmooth regularization function. Each agent has access only to its individual objective function, and performs updates using local data. A network-wide solution is achieved asymptotically via an exchange of information among neighboring agents, which serves to effectively integrate global information. 

Distributed first order methods \cite{1st2009}-\nocite{1st2015}\nocite{1st2018}\cite{usman2} constitute a popular class of algorithms for solving (\ref{prob1}) due to their simplicity and low computational requirements. Agents perform updates using (sub) gradients pertaining to local objectives along with averaged values from neighbors. The method in \cite{1st2009} can be viewed as solving a penalized version of (\ref{prob1}) and therefore converges to a biased solution. To converge to the exact solution, diminishing step sizes have to be used and consequently converges slowly. This is remedied by the method in \cite{1st2018}, where agents utilize iterates and gradient histories. By using tracking techniques, \cite{1st2018} establish $\mathcal{O}(1/t)$ convergence rate under general convex objectives and linear convergence rate under strongly convex objectives. Despite its economical costs, first order methods are well-known to suffer from slow convergence speed, especially for ill-conditioned large scale problems. 

Second order methods attain acceleration by invoking information pertaining to the (exact or approximate) Hessian matrix of the objective. Recent work on (quasi) Newton methods \cite{nn2017}-\nocite{nn2020}\nocite{dnadmm}\cite{BFGS2021} incorporate objective curvature information to compute more informed updates. Similar to its first order counterpart, Network Newton \cite{nn2017} and its asynchronous extension \cite{nn2020} solve a penalized version of (\ref{prob1}), where Newton steps are approximated via a truncated Taylor expansion. Similar techniques were employed in \cite{dnadmm} for the primal-dual formulation of the problem. The ability to reach a high-accuracy solution in significantly fewer iterations makes second order algorithms a desirable candidate when abundant computational resources are available. However, each update step of Newton's method requires to solve a linear system at the cost of $\mathcal{O}(d^3)$ computational costs for general objectives. 

Methods that allow agents to choose between first and second order solver are rather scarce in the literature. This consideration becomes significant for large scale heterogeneous networks where machines have varying capacities in terms of computation, communication, storage, and battery levels. In this aspect, a \emph{hybrid} scheme serves to effectively avoid placing heavy computational burden on less capable machines thus better utilizing the network-wide resources. An emerging application where such challenges become more pronounced is Federated Learning \cite{FL2017}, where massive number of agents collaborate and aim to obtain a model with improved performances by utilizing data distributed across the network. We note a related work \cite{fedhybrid2021} where authors proposed a hybrid protocol for the client-server setting where clients may choose to perform either gradient or Newton updates depending on their local environments. Nevertheless, \cite{fedhybrid2021} is limited to the client-server setting, and further requires synchronous computation and communication: this is unrealistic in large-scale networks. In contradistinction, the proposed algorithm can be deployed in general multi-agent settings, (of which client-server is a special case) and achieves a totally asynchronous implementation.  

\noindent \textbf{Contributions}:
\begin{itemize}
    \item We propose \underline{H}ybr\underline{i}d \underline{P}rimal-dual \underline{P}roximal \underline{O}ptimization (HIPPO), a framework that encases heterogeneous multi-agent systems by allowing agents to select between gradient and or Newton updates. 
    \item HIPPO admits a \textit{totally asynchronous} implementation that allows a general participation scheme of agents at its iteration. This removes the need for a central clock and network-wide synchronization. 
    \item We establish linear convergence rate in expectation under strongly convex local objectives with Lipschitz continuous gradients, for arbitrarily chosen hybrid updating schemes. Numerical experiments on real-life datasets attest the efficacy of the proposed algorithm. 
\end{itemize}
\noindent \textbf{Notation}: Vectors $x\in\mathbb{R}^d$ are denoted using lower case letters and matrices $A\in\mathbb{R}^{n\times m}$ are denoted using upper case letters, while transpose is denoted as $A^\top$. Row stacking and column stacking is expressed as $[A;B]$ and $[A,B]$, respectively (provided dimensions match). We use superscript to denote iterates of the algorithm, while subscripts are reserved for vector components. For example, $x_i^t$ denotes the vector held by the $i$-th agent at time step $t$. When a norm specification is not provided, $\norm{x}$ and $\norm{A}$ denote the Euclidean norm of the vector and the induced norm of the matrix, respectively. For a positive definite matrix $P\succ 0$, we represent the associated quadratic form as $\norm{x}_P:=\sqrt{x^\top P x}$. For matrices of arbitrary dimensions, we denote their Kronecker product as $A\otimes B$. Last, we use the shorthand notation $[m]:=\{1,\dots,m\}$.  

\section{PRELIMINARIES}
\subsection{Proximal operator} 
The proximal operator \cite{convex_analysis} $\textbf{prox}_{g/\mu}(\cdot):\mathbb{R}^d\to \mathbb{R}$ associated with a function $g(\cdot):\mathbb{R}^d\to \mathbb{R}$ and a constant $\mu>0$ is defined as: 
\begin{align}
    \textbf{prox}_{g/\mu}(x) &=\underset{v\in\mathbb{R}^d}{\text{argmin}}\left\{g(v)+\tfrac{\mu}{2}\norm{v-x}^2\right\}. \label{proximal}
\end{align}
When $g(\cdot)$ is convex, the above mapping is single valued. Proximal algorithms \cite{proximal_algo} were developed for convex composite problems, where the objective is the sum of two convex functions: one smooth and one possibly non-smooth but proximable (i.e., with efficiently computable proximal mapping). Examples where (\ref{proximal}) is easy to evaluate are ubiquitous in applications, e.g., the $\ell_1-$norm or the  indicator functions for hyperplane, halfspace, simplex, orthant, etc.  
The vast majority of these methods are first order, i.e., they use the gradient pertaining to the smooth term along with the proximal operator of the non-smooth part to iteratively solve the problem. Using second order information, in this framework, typically poses challenges in ensuring convergence guarantees while maintaining affordable computational costs. This is due to the fact that the norm function in (\ref{proximal}) has to be scaled according to the Hessian matrix \cite{proximal_newton}. Such scaling may render the evaluation of the associated proximal mappings expensive, even when its first-order counterpart admits a closed form solution.

\subsection{Problem formulation}
The network topology is captured by an undirected connected graph: $\mathcal{G}:= (\mathcal{V},\mathcal{E})$, where $\mathcal{V}:=[m]$ is the index set and $\mathcal{E}\subseteq (\mathcal{V},\mathcal{V})$ is the edge set. Edge $(i,j)\in \mathcal{E}$ if and only if a communication channel exists between the $i$-th agent and $j$-th agent. We do not consider self-loops, i.e., $(i,i)\notin \mathcal{E}$ for all $i\in[m]$. The number of edges is represented as $\abs{\mathcal{E}}=n$, and the set of neighbors of the $i$-th agent is defined by $\mathcal{N}_i:=\{j\in [m], (i,j)\in \mathcal{E}\}$. Using the above definitions, problem (\ref{prob1}) can be recast as a consensus optimization problem by introducing local decision vectors $x_i$ to agents, and consensus vectors $z_{ij}$ to edges of the graph. In specific: 
\begin{equation}\label{prob2}
\begin{aligned}
    \underset{x_i,\theta,z_{ij}\in\mathbb{R}^d}{\mathrm{minimize}}\,\, \bigg\{\sum_{i=1}^m f_i(x_i) & + g(\theta) \bigg\} \\
    \text{subject to}\,\,x_i= z_{ij}=  x_j,\,\,&\text{for  $i\in [m], j\in \mathcal{N}_i$.,}\\
    x_l =\theta \,\,\text{for}\,\,&\text{some $l\in[m]$.}
\end{aligned}
\end{equation}
Note that we have also separated the argument of the smooth part and the nonsmooth part of the objective by letting $x_l=\theta$ for an arbitrary agent. The introduction of intermediate variables $\{z_{ij}\}$ deviates from the classical setting of consensus optimization and is crucial in obtaining block-diagonal Hessian structures that are needed for our hybrid scheme (See Section \ref{hybrid section}). The above formulation can be compactly expressed by defining the source and the destination matrices $\widehat{A}_s,\widehat{A}_d$ as follows. We first assign an arbitrary order to the edge set and denote $\mathcal{E}_k=(i,j), k\in [n],$ if the $k$-th edge connects the $i$-th and the $j$-th agent. Matrices $\widehat{A}_s,\widehat{A}_d\in\mathbb{R}^{n\times m}$ are constructed by letting $[\widehat{A}_s]_{ki}=[\widehat{A}_d]_{kj}=1$ (with all other entries in the $k-$th row being zero). The network topology can be captured by the signed and unsigned graph incidence matrices defined respectively as: $\widehat{E}_s=\widehat{A}_s-\widehat{A}_d$ and $\widehat{E}_u=\widehat{A}_s+\widehat{A}_d$. An immediate consequence is that problem (\ref{prob2}) is equivalent to: 
\begin{equation}\label{compact}
\begin{aligned}
     \underset{x\in \mathbb{R}^{md},z\in\mathbb{R}^{nd},\theta\in\mathbb{R}^d}{\mathrm{minimize}}\,\, \bigg\{F(x)&+g(\theta)\bigg\} \\
     \text{subject to}\,\, \begin{bmatrix}A_s \\A_d \end{bmatrix}x= Ax &=Bz = \begin{bmatrix}I_{nd}\\ I_{nd}  \end{bmatrix}z,\\
     S^\top x &= \theta.
    \end{aligned}
\end{equation}
where we stack $x_i$ into a single column vector $x\in\mathbb{R}^{md}:=[x_1^\top, \dots, x_m^\top]^\top$ and similarly for $z\in\mathbb{R}^{nd}$. The matrix $S:=(s_l\otimes I_d)\in\mathbb{R}^{md\times d}$, where $s_l\in\mathbb{R}^d$ is a zero vector except for the $l$-th entry being one, serves to select the $l$-th agent. We denote $A_s:=\widehat{A}_s\otimes I_d,A_d:=\widehat{A}_d\otimes I_d$ as the Kronecker product between $\widehat{A}_s,\widehat{A}_d$ and the identity matrix of dimension $d$, respectively. Matrices $A,B$ in (\ref{compact}) are formed by stacking the corresponding matrices. 

\subsection{Primal-dual algorithms}
Primal-dual optimization algorithms date back to the work \cite{earlyadmm}. By introducing the (augmented) Lagrangian, primal-dual algorithms seek to solve the primal and the dual problem simultaneously, through the use of dual variables. The proliferation of primal-dual methods \cite{bert} is attributed to their wide applicability to distributed optimization. In the following of this subsection, we review one instance of primal-dual algorithms, termed Alternating Direction Method of Multipliers (ADMM) \cite{boydadmm}, which serves as the foundation of our subsequent development. To solve (\ref{compact}) using ADMM, we define the associated augmented Lagrangian as follows: 
\begin{gather}
    \mathcal{L}(x,z,\theta;y,\lambda) = F(x)+g(\theta)+y^\top (Ax-Bz)\nonumber\\
    +\lambda^\top (S^\top x-\theta) +\tfrac{\mz}{2}\norm{Ax-Bz}^2+\tfrac{\mt}{2}\norm{S^\top x -\theta}^2,  \label{AL}
\end{gather}
where $y\in\mathbb{R}^{2nd},\lambda\in\mathbb{R}^d$ are the dual vectors and $\mz,\mt>0$ is are scalars that controls the quadratic penalty. ADMM solves (\ref{compact}) by sequentially minimizing the augmented Lagrangian with respect to the primal variables $(x,z,\theta)$, and then performs gradient ascent on the dual variable. 
\begin{subequations}\label{3blockadmm}
\begin{align}
    x^{t+1} &= \underset{x}{\mathrm{argmin}}\,\,\mathcal{L}(x,z^t,\theta^t;y^t,\lambda^t), \label{admm1}\\
    z^{t+1} &=\underset{z}{\mathrm{argmin}}\,\,\mathcal{L}(x^{t+1},z,\theta^t;y^t,\lambda^t), \label{admm2}\\
    \theta^{t+1} &=\underset{\theta}{\mathrm{argmin}}\,\,\mathcal{L}(x^{t+1},z^{t+1},\theta;y^t,\lambda^t)\label{admm3}\\
    y^{t+1} &= y^{t}+\mz (Ax^{t+1}-Bz^{t+1}). \label{admm4}\\
    \lambda^{t+1} &= \lambda^t+\mt(S^\top x^{t+1}-\theta^{t+1})\label{admm5}
\end{align}
\end{subequations}
As a result of introducing the consensus variables $\{z_{ij}\}$, we obtain an instance of 3-block ADMM in (\ref{3blockadmm}), which impose more strict assumptions than 2-block ADMM for convergence\cite{3block}. The introduction of the consensus variables $\{z_{ij}\}$ decouples $\{x_i\}$ and renders the step (\ref{admm1}) distributedly computable \cite{wei2014}, i.e., agents can perform the step (\ref{admm1}) with only local information and information from their immediate neighbors. However, multiple inner loops are required to solve (\ref{admm1}) with general objective functions. To alleviate the computational burden, several approximation schemes \cite{Ling2015}-\nocite{mokhtari2016}\cite{lbfgs2022} exist wherein instead of aiming for an exact solution, agents perform one or more updating steps using gradient or (quasi) Newton information. Building on these results, we propose a general framework that accommodates system heterogeneity by allowing agents to select their updates and further compute updates in an asynchronous fashion. 

\section{ALGORITHMIC DEVELOPMENT}\label{section3}
Abbreviating $\mathcal{L}^t\equiv \mathcal{L}(x^t,z^t,\theta^t;y^t,\lambda^t)$, we define the approximated augmented Lagrangian as follows, 
\begin{gather}
    \widehat{\mathcal{L}}(x,z^t,\theta^t;y^t,\lambda^t)=\mathcal{L}^t+(x-x^t)^\top \nabla_x \mathcal{L}^t
    +\tfrac{1}{2}\norm{x-x^t}^2_{H^t}, \label{AL_model}
\end{gather}
where $H^t$ is a matrix (to be defined subsequently) that encapsulates our hybrid updating schemes. By replacing the augmented Lagrangian (\ref{AL}) with its approximation in (\ref{AL_model}), we obtain the following closed form solution for step (\ref{admm1}):
\begin{gather}
    x^{t+1} = x^t- (H^t)^{-1} \nabla_x \mathcal{L}(x^t,z^t;y^t). \label{x_primal}
\end{gather}
Moreover, step (\ref{admm2}) corresponds to: 
\begin{gather}
    B^\top y^t+\mu B^\top (A x^{t+1}-Bz^{t+1})= 0 .\label{z_update}
\end{gather}
By completion of squares, the update of $\theta^{t+1}$ can be cast as the following proximal mapping: 
\begin{align}
    \theta^{t+1} = \textbf{prox}_{g/\mt}(S^\top x^{t+1}+\tfrac{1}{\mt}\lambda^t). \label{theta_update}
\end{align}
Dual variables $(y,\lambda)$ are updated in verbatim as (\ref{admm4})-(\ref{admm5}). The updates (\ref{x_primal})-(\ref{theta_update}), and (\ref{admm4})-(\ref{admm5}) serve as the basis for developing our proposed algorithm. An equivalent while more efficient implementation can be achieved, which we formalize in the following lemma.     

\begin{lemma}\label{lemma1} Consider the iterates generated by (\ref{x_primal})-(\ref{theta_update}), and (\ref{admm4})-(\ref{admm5}). Let the dual variable $y^t\in \mathbb{R}^{2nd}$ be expressed as $y^t=[\alpha^t;\beta^t],\alpha^t,\beta^t\in\mathbb{R}^{nd}$. If $y^0$ and $z^0$ are initialized so that $\alpha^0=-\beta^0$ and $z^0=\tfrac{1}{2}E_ux^0$, then $\alpha^t=-\beta^t$ and $z^t=\tfrac{1}{2}E_u x^t$ for all $t\geq 0$. Moreover, by defining $\phi^t=E_s^\top \alpha^t\in \mathbb{R}^{md}$, we obtain the following equivalent updates:
\end{lemma}
\begin{subequations}\label{update}
\begin{align}
    x^{t+1} &= x^t- (H^t)^{-1}\bigg(\nabla F(x^t)+\phi^t+\tfrac{\mz}{2}E_s^\top E_s x^t+S\Big[\lambda^t\nonumber\\
    &\,\,+\mt(S^\top x^t-\theta^t)\Big]\bigg),\label{lemma1_x}\\
    \theta^{t+1} &= \textbf{prox}_{g/\mu} (S^\top x^{t+1}+\tfrac{1}{\mt}\lambda^t), \label{lemma1_theta}\\
    \phi^{t+1} & = \phi^t +\tfrac{\mu}{2}E_s^\top E_s x^{t+1},\label{lemma1_dual}\\
    \lambda^{t+1} &= \lambda^t+\mt(S^\top x^{t+1}-\theta^{t+1}). \label{lemma1_lambda}
\end{align}
\end{subequations}
\textit{Proof}: See Appendix. 
\begin{remark}
The above result is similarly derived in \cite{wei2014} when the subproblem is solved exactly. Note updates (\ref{update}) are equivalent to updates (\ref{x_primal})-(\ref{theta_update}), and (\ref{admm4})-(\ref{admm5}) under appropriate initialization, irrespective of choices of $H^t$. However, (\ref{update}) falls into the category of 2-block ADMM through elimination of the $z$-update ($z^t=\tfrac{1}{2}E_u x^t$). This gives more freedom for choosing hyperparameters $\mz,\mt$. On the other hand, the use of the consensus variables $\{z_{ij}\}$ in the formulation (\ref{compact}) is necessary, as we demonstrate next.\end{remark}

\begin{algorithm}[t]
    \caption{HIPPO} 
    Zero initialization for all variables. Hyperparameters $\mz,\mt$ and $\Delta_{ii},i\in[m]$. 
    \begin{algorithmic}[1]
    \For{$t=0,1,2,\ldots$}
    \For{active agents}
    \State \label{retrieve}Retrieve $x_j^t,j\in\mathcal{N}_i$, from the local buffer.
    \label{comp1}\If{agent $i\in \mathcal{U}^t_{\mathrm{Gradient}}$}
        \State $J_{ii}^t = 0$
        \ElsIf{$i\in \mathcal{U}^t_{\mathrm{Newton}}$}
        \State $J_{ii}^t = \nabla^2 f_i(x_i^t)$
    \EndIf
    \label{comp2}\State Compute $u_i^t$ as in (\ref{u_direction}).
    \State \label{primal_algo}$x^{t+1}_i=x^t_i-u^t_i$.
    \State \label{broadcast_algo}Broadcast $x_i^{t+1}$.
    \State \label{dual_algo}$\phi_i^{t+1} = \phi_i^t+\tfrac{\mu}{2}\sum_{j\in \mathcal{N}_i}(x_i^{t+1}-x_j^{t+1})$.
    \If{$i=l$}
		\State \label{proximal_start} $\theta^{t+1}=\textbf{prox}_{g/\mt}(x_l^{t+1}+\tfrac{1}{\mt}\lambda^t)$
		\State \label{proximal_end}$\lambda^{t+1}=\lambda^t+\mt(x^{t+1}_l-\theta^{t+1})$
	\EndIf 
    \EndFor
    \EndFor
\end{algorithmic}
\end{algorithm}

\subsection{Hybrid updating}\label{hybrid section}
We construct $H^t$ as follows: 
\begin{gather}
    H^t = J^t + \mz D+\mt SS^\top+ \Delta, \label{H_mtx}
\end{gather}
where $J^t\in \mathbb{R}^{md\times md}$ is a block diagonal matrix, whose sub-blocks dictate the updating scheme of the corresponding agent, $D$ is a diagonal matrix with entries equal to the degree of the corresponding agent, and $\Delta$ is a constant matrix, whose sub-blocks serve to specify step sizes and provide numerical robustness. In specific, if the $i$-th agent opts to use Newton updates, then $J^t_{ii}=\nabla^2 f_i(x_i^t)\in \mathbb{R}^{d\times d}$; if the $i$-th agent opts to use first-order information, then $J^t_{ii}=0$ and the step size of the gradient descent is determined by specifying $\Delta_{ii}$. Note that when $J^t=\nabla^2 F(x^t)$, i.e., all agents perform Newton updates, we obtain the Hessian of the augmented Lagrangian (\ref{AL}) plus the diagonal step size matrix $\Delta$. Irrespective of the choice of $J^t_{ii}$, $H^t$ in (\ref{H_mtx}) is always block diagonal. This important property translates to the fact that computing the update step $u^t=(H^t)^{-1} \nabla_x \mathcal{L}^t$ can be carried out \textit{locally} by agents, without additional communication rounds once $\nabla_x \mathcal{L}^t$ is available. In specific, the $i$-th agent computes the update direction $u_i^t$ by solving the following linear system: 
\begin{gather}
    (J^t_{ii}+\mz \abs{\mathcal{N}}_i+\delta_{il}\mt+\Delta_{ii}) u_i^t = \nabla f_i(x_i^t)+\phi_i^t\nonumber\\
    +\tfrac{\mz}{2} \sum_{j\in \mathcal{N}_i} (x_i^t-x_j^t)+\delta_{il}(\lambda^t+\mt(x^t_l-\theta^t)),\label{u_direction}
\end{gather}
where $\delta_{il}=1$ if and only if $i=l$. The block-diagonal structure of $H^t$ in (\ref{H_mtx}) is a consequence of the introduction of $\{z_{ij}\}$. If instead a direct consensus constraint is imposed as $x_i=x_j$ in (\ref{prob2}), it is not hard to verify that $H^t$ in (\ref{H_mtx}) would be of the following form: 
\begin{gather*}
    J^t+\mz L_s+\mt SS^\top +\Delta,
\end{gather*}
where $L_s=E_s^\top E_s$ is the signed graph Laplacian matrix, whose $ij$-th block is nonzero if $(i,j)\in \mathcal{E}$. The presence of $L_s$ couples agents with their neighbors. In such cases, solving for the update direction $u_i^t$ requires either collecting all information from neighbors beforehand, or invoking a truncated Taylor expansion of the Hessian inverse as in \cite{dnadmm}. Such schemes not only induce additional communication rounds within agents, but are actually infeasible in the hybrid scheme: if neighboring agents opt to use different updates, the required information to compute approximated Newton steps is not available. The proposed algorithm is described in \textbf{Algorithm 1}.

We equip each agent with a local buffer so that agents can receive and store information from neighbors even if the agent is not active. We assume each agent is active according to local Poisson clocks \cite{nn2020}. For a round that agent $i$ participates, i.e., performs an update, it retrieves the most recent copy of $x_j^t,j\in\mathcal{N}_i$, from its local buffer as in line \ref{retrieve} of \textbf{Algorithm 1}. Active agents are divided into two groups: $\mathcal{U}^t_{\mathrm{Gradient}}$ and $\mathcal{U}^t_{\mathrm{Newton}}$, where $\mathcal{U}^t_{\mathrm{Gradient}}\cup\, \mathcal{U}^t_{\mathrm{Newton}}\subseteq \mathcal{V} $. Depending on whether an agent opts to use gradient or Newton update, the update direction is computed using $J^t_{ii}=\nabla^2 f(x_i^t)$ or $J^t_{ii}=0$ in (\ref{u_direction}) accordingly (lines \ref{comp1}-\ref{comp2}). Note that each agent can choose different computing schemes at each time step. Once the update direction $u_i^t$ is obtained, the primal update is carried out, followed by broadcasting $x_i^{t+1}$ to neighbors in lines \ref{primal_algo}-\ref{broadcast_algo}. For Poisson processes, it is natural to assume only one clock ticks at a time and therefore $x_j^{t+1}=x_j^t$. For other cases, agents can check their local buffers before proceeding to the dual update in the line \ref{dual_algo}. The $l$-th agent additionally update the $(\theta,\lambda)$ pair pertaining to the regularization function.

\subsection{Operator description of the algorithm}\label{operator_subsection}

To succinctly describe the proposed asynchronous algorithm for the sake of analysis, we introduce the following vector concatenation $v:= [x; \phi;\theta;\lambda]\in\mathbb{R}^{(2m+2)d}$. We define an operator $T:\mathbb{R}^{(2m+2)d}\to \mathbb{R}^{(2m+2)d}$ as follows:
\begin{align}
    v^{t+1} = Tv^t, \label{sync_operator}
\end{align}
which maps $[x^t;\phi^t;\theta^t;\lambda^t]\in\mathbb{R}^{(2m+2)d}$ to $[x^{t+1};\phi^{t+1};\theta^{t+1};\lambda^{t+1}]$ according to (\ref{update}). We further define the following activation matrix: 
\begin{align}
    \Omega^{t} = 
    \begin{bmatrix}
    X^{t} & 0         & 0           & 0\\
    0       & X^{t}   & 0           & 0\\
    0       & 0         &X^{t}_{ll} & 0\\
    0       & 0         & 0           & X^{t}_{ll} 
    \end{bmatrix}, \label{omega_mtx}
\end{align}
where $X^{t}\in \mathbb{R}^{md\times md}$ is a random diagonal matrix with sub-blocks $X^{t}_{ii}\in\mathbb{R}^{d\times d},i\in[m]$, taking values $I_d$ or the zero matrix. Using (\ref{sync_operator})--(\ref{omega_mtx}), we can express the asynchronous updates as: 
\begin{align}
    v^{t+1}= v^t+\Omega^{t+1}(Tv^t-v^t). \label{asy_update}
\end{align}
This update rule simply prescribes that the $i$-th agent participates at the $t$-th iteration if and only if $X_{ii}^{t+1}=I_d$. Moreover, the primal-dual pair $(x_i^{t+1},\phi_i^{t+1})$ is updated simultaneously when the corresponding agent is active, by construction. Therefore, by varying $\mathbb{E}^t\left[X_{ii}^{t+1}\right]=\mathbb{P}(X_{ii}^{t+1}=I_d\vert \mathcal{F}^t)$, where $\mathcal{F}^t$ is the filtration generated by $X^1,X^2,\dots,X^t$ and $\mathbb{E}^t[\cdot]:=\mathbb{E}[\cdot\vert \mathcal{F}^t]$, one obtain a wide range of activation schemes. 
\begin{itemize}
    \item \textit{Synchronous updates}: $X_{ii}^t=I_d$ for $i\in[m]$ and $t\geq 1$. In this case, updates (\ref{asy_update}) boils down to (\ref{sync_operator}). 
    \item \textit{Stochastic coordinate updates}: By letting the random matrix $X^{t}$ take values in:\{ $\textbf{Blkdiag}(I_d,0,\dots,0)$, $\textbf{Blkdiag}(0,I_d,0,\dots,0)$, $\dots$, $\textbf{Blkdiag}(0,0,\dots,I_d)$\}, we obtain a randomized algorithm where, at each iteration, a single agent performs gradient or Newton update, depending on $i\in\mathcal{U}^t_{\mathrm{Gradient}}$ or $i\in\mathcal{U}^t_{\mathrm{Newton}}$. 
    \item \textit{Multi-agent updates}: By setting $\mathbb{E}^t[X^{t+1}_{ii}]=p_i^t$, we obtain a general activation scheme where the $i$-th agent performs updates at $t$-th iteration with probability $p_i^t$. 
\end{itemize}

\section{ANALYSIS}
In this section, we present a unified analysis that establishes linear convergence rate for the proposed algorithm under strongly convex local objectives, with arbitrary updating schemes chosen by agents. We assume that the initialization requirement in Lemma 1 is satisfied (zero initialization for all variables suffices). Recall the definition $\phi= E_s^\top \alpha$ in Lemma 1 and the concatenation $v=[x;\phi;\theta;\lambda]\in\mathbb{R}^{(2m+2)d}$ in Section \ref{operator_subsection}. Using the tuple $(x,\phi,\theta,\lambda)$ facilitates the implementation of the algorithm while $(x,z,\alpha,\theta,\lambda)$ is most suitable for our analysis. We note that their equivalence is established in Lemma 1. In the following, we denote $v_\alpha=[x;z;\alpha;\theta;\lambda]\in\mathbb{R}^{(m+2n+2)d}$. The following assumptions are made throughout.  

\begin{assumption}\label{a1}
The local cost functions $f_i(\cdot):\mathbb{R}^d\to\mathbb{R}$ are twice continuously differentiable, $m_f-$strongly convex, and have $M_f-$Lipschitz continuous gradients, i.e., $\forall\,i\in[m], x\in\mathbb{R}^d$, the following holds:
\end{assumption}
\begin{align}
    m_f I \preceq \nabla^2 f_i(x)\preceq M_f I. \label{assump1}
\end{align}
\begin{assumption}
The function $g(\cdot):\mathbb{R}^d\to \mathbb{R}$ is proper, convex, and closed, i.e., $\forall x,y\in\mathbb{R}^d$, the following holds:
\end{assumption}
\begin{align}
    (\partial g(x)-\partial g(y))^\top (x-y)\geq 0,
\end{align}
where $\partial g(\cdot)$ denotes the subdifferential set.

\begin{assumption}\label{a3}
The Hessian of the local cost function $f_i(\cdot)$ is Lipschitz continuous with constant $L_f$, i.e., $\forall\,i\in[m], x,y\in\mathbb{R}^d$, the following holds:
\end{assumption}
\begin{align}
    \norm{\nabla^2 f_i(x)-\nabla^2 f_i(y)}\leq L_f\norm{x-y}.
\end{align}
Since $F(x)=\sum_{i=1}^m f_i(x_i)$, the bound in (\ref{assump1}) holds for $\nabla^2 F(x)$ as well, i.e., $m_f I\preceq  \nabla^2 F(x)\preceq M_f I$. The Lipschitz continuity assumption is standard in analyzing both first and second order methods while Assumption \ref{a3} is only needed for analyzing second order methods \cite{convex_opt}. We proceed by stating a lemma that that will be handy in quantifying the error introduced when we replace the exact optimization step (\ref{admm1}) with the one-step update (\ref{lemma1_x}).

\hide{We proceed to state the KKT condition for the tuple $(x,z,\alpha)$ in the following. 
\begin{lemma}
Recall the identity $E_s=A_s-A_d$ and $E_u=A_s+A_d$. The tuple $(x^\star,z^\star,\theta^\star,\alpha^\star,\lambda^\star)$ solves (\ref{compact}), and equivalently (\ref{prob1}) if and only if the following holds: 
\end{lemma}
\begin{align}
    \nabla F(x^\star)+E_s^\top \alpha^\star+S \lambda^\star &= 0, \tag*{KKTa}\label{KKTa}\\
    \partial g(\theta^\star)-\lambda^\star &\ni 0 ,\tag*{KKTb}\label{KKTb}\\
    E_s x^\star& =0,\tag*{KKTc}\label{KKTc}\\
    E_u x^\star& = 2z^\star ,\tag*{KKTd}\label{KKTd}\\ 
    S^\top x^\star&=\theta^\star.\tag*{KKTe}\label{KKTe}
\end{align}
\textit{Proof}: See Appendix. }

\begin{lemma}\label{lemma2}
Consider the iterates generated by (\ref{update}) and equivalently (\ref{sync_operator}). Recall the definition $\phi=E_s^\top \alpha^t$ and the identity $z^t=\tfrac{1}{2}E_u x^t$. The following holds:
\end{lemma}
\scalemath{0.95}{
\begin{aligned}
    e^t+\nabla F(x^{t+1})-\nabla F(x^\star)+\Delta (x^{t+1}-x^t)+E_s^\top (\alpha^{t+1}-\alpha^\star)\nonumber\\
    +\mz E_u^\top(z^{t+1}-z^t)
   +S\left( \lambda^{t+1}-\lambda^\star+\mt\left( \theta^{t+1}-\theta^t\right)\right) = 0,\nonumber
\end{aligned}
}
where the error term $e^t$ is defined as: 
\begin{gather}
        e^t := \nabla F(x^t)-\nabla F(x^{t+1})
        +J^t(x^{t+1}-x^t). \label{error}
\end{gather}
\textit{Proof}: See Appendix. 

The following lemma establishes an upper bound for $\norm{e^t}$ that is crucial to establish the convergence of the algorithm. 
\begin{lemma}\label{error_bound_lemma}
The error $e^t$ in (\ref{error}) is upper bounded as: $\norm{e^t}^2\leq \norm{x^{t+1}-x^t}^2_{(\Pi^t)^2}$, where $\Pi^t\in\mathbb{R}^{md\times md}$ is a diagonal matrix with entries equal to: 
\end{lemma}
\begin{align}
    \Pi_{ii}^t &= M_f \cdot I,\,\,\text{for}\,\,i\in \mathcal{U}^t_{\mathrm{Gradient}}, \label{lemma3_eqn}\\
    \Pi_{ii}^t &= \min\left\{2M_f, \tfrac{L_f}{2}\norm{x_i^{t+1}-x_i^t}\right\}\cdot I,\,\,\text{for}\,\,i\in\mathcal{U}^t_{\mathrm{Newton}}.\nonumber
\end{align}

Lemma \ref{error_bound_lemma} reveals the difference between using Newton and gradient updates, i.e., the progression of the error in terms of $\norm{x^{t+1}-x^t}$. Before we establish the convergence rate of the algorithm, we introduce the following activation matrix for the tuple $(x,z,\alpha,\theta,\lambda)$:
\begin{gather}
    \Omega_\alpha^{t}: =\begin{bmatrix}
    X^{t}   & 0   &        0 &       0 & 0\\
    0       & Y^t &        0 &       0 & 0 \\
    0       & 0   &      Y^t &       0 & 0 \\
    0       &   0 &        0 &X^t_{ll} & 0 \\
    0       &   0 &        0 & 0       &X^t_{ll}
    \end{bmatrix}.\label{activation_mtx_2}
\end{gather}
where we have additionally introduced a diagonal random matrix $Y^t\in\mathbb{R}^{nd\times nd}$. The activation matrix $\Omega^t_\alpha$ in (\ref{activation_mtx_2}) differs from $\Omega^t$ in (\ref{omega_mtx}) as we allow $(z,\alpha)$ to be updated independently, as captured by the activation matrix $Y^t$. Similar to (\ref{asy_update}), we compactly write the iterates as: 
\begin{align}
    v^{t+1}_\alpha = v^t_\alpha+\Omega^{t+1}_\alpha (Tv^t_\alpha-v^t_\alpha). \label{asy_update_2}
\end{align}

To compare (\ref{asy_update}) (\textbf{Algorithm 1}) and (\ref{asy_update_2}), we assign $(z,\alpha)$ to edges and $(x,\theta,\lambda)$ to agents (note that $(\theta,\lambda)$ is only held by the $l$-th agent). We emphasize that this is done \emph{only} for analysis, and \textbf{Algorithm} 1 does not require storing/updating edge variables. In fact, the activation scheme in (\ref{asy_update}) is a special case of (\ref{asy_update_2}), and has a simple realization: whenever an agent is active, all incident edges become active and update their variables. We first characterize the condition that ensures the convergence of (\ref{asy_update_2}) in Theorem 1 and then show that (\ref{asy_update}) satisfies this condition. 

\begin{theorem}\label{theorem1}
Consider the iterates generated by (\ref{asy_update_2}), and let assumptions 1-3 hold. We denote the maximum and minimum eigenvalue of $L_u=E_u^\top E_u$ as $\sigma^{L_u}_\mathrm{max}$ and $\sigma^{L_u}_\mathrm{min}$ respectively. Under any activation schemes such that $\mathbb{E}^t[\Omega^{t+1}_\alpha]=\Omega_\alpha \succ 0$, set $\mz=2\mt$, $\Delta_{ii}=\epsilon$ for $i\in[m]$. Then the following holds:
\end{theorem}
\begin{gather*}
    \mathbb{E}^t\left[\norm{\va^{t+1}-\vas}^2_{\mathcal{H}\Omega_\alpha^{-1}}\right]\leq \left(1-\tfrac{p^{\mathrm{min}}\eta}{1+\eta}\right)\norm{\va^{t}-\vas}^2_{\mathcal{H}\Omega_\alpha^{-1}},
\end{gather*}
where for $i\in[m],k\in[n]$, $\mathbb{E}^t[X^{t+1}_{ii}]=p^X_i,\mathbb{E}^t[Y^{t+1}_{kk}]=p^Y_k$, $p^{\mathrm{min}}:=\underset{i\in[m],k\in[n]}{\min}\{p_i^X,p_k^Y\}$, $\eta$ is given by:
\begin{gather}
    \eta=\min \bigg\{\tfrac{2m_fM_f}{m_f+M_f}\tfrac{1}{\epsilon+\mt(\sigma^{L_u}_{\mathrm{max}}+2)},\tfrac{1}{2},\tfrac{2}{5}\tfrac{\mt\sigma^+_{\mathrm{min}}}{m_f+M_f},\nonumber\\
    \tfrac{\mt\sigma^+_{\mathrm{min}}}{5\epsilon},\tfrac{\sigma^+_{\mathrm{min}}}{5\max\{1,\sigma^{L_u}_{\max}\}}\bigg\}\label{eta},
\end{gather}
and the scaling matrix $\mathcal{H}=\textbf{diag}[\epsilon,2\mz,\tfrac{2}{\mz},\mt,\tfrac{1}{\mt}]$.\\
\textit{Proof}: See Appendix. 

Note that the only requirement on the activation scheme is $\mathbb{E}^t[\Omega^{t+1}_\alpha]\succ 0$, i.e., each agent and each edge has a nonzero probability of being active at each round. We proceed to show that by specifying $X^{t+1}$ alone (sampling agents and activating incident edges) in (\ref{asy_update}), the corresponding tuple $(x,z,\alpha,\theta,\lambda)$ converges linearly in expectation.
\begin{corollary}
Consider the iterates generated by (\ref{asy_update}). Under the same setting as in Theorem \ref{theorem1}, for arbitrary activation scheme such that $\mathbb{E}^t[X^{t+1}]\succ 0$, the corresponding tuple $(x,z,\theta;\alpha,\lambda)$ converges linearly in expectation. 
\end{corollary}
\textit{Proof}: See Appendix. 

\begin{remark}
Note that Theorem 1 and Corollary 1 does not impose any conditions on the proportion of agents who perform gradient or Newton updates. This is due to the fact that we can upper bound the matrix $\Pi^t$ in Lemma  \ref{error_bound_lemma} as $\Pi^t\preceq 2M_f I$. To quantify the impact on the proportion of agents performing Newton updates, we consider a special case of activation scheme in (\ref{asy_update_2}): we activate one variable uniformly at random with probability $p=\tfrac{1}{m+n}$, so that $\mathbb{E}^t[\Omega^{t+1}_\alpha]=pI$. Moreover, we denote the percentage of agents who perform Newton updates as $q\in [0,1]$. Then the following holds: 
\end{remark}
\begin{gather*}
    \mathbb{E}^t\left[\norm{e^{t+1}_v}^2\right]
    \leq (1-q)\left(1-\tfrac{p\eta_{\mathrm{g}}}{1+\eta_{\mathrm{g}}}\right)\norm{e^t_v}^2\\
    + q\left(1-\tfrac{p\eta_{\mathrm{n}}}{1+\eta_{\mathrm{n}}}\right)\norm{e^t_v}^2,
\end{gather*}
where $\norm{e^t_v}^2=\norm{v^t_\alpha-v^\star_\alpha}^2$ and $\eta_{\mathrm{g}},\eta_{\mathrm{n}}$ denotes the convergence coefficient of gradient and Newton updates respectively. Recall Lemma \ref{error_bound_lemma}. Note that as $t\to\infty$, $\Pi^t_{ii}\to 0$ for $i\in \mathcal{U}_{\mathrm{Newton}}$. This allow us to choose $\Delta_{ii}=0$ for $i\in \mathcal{U}_{\mathrm{Newton}}$. By inspecting (\ref{hybrid_eta}), we conclude that $\eta_{\mathrm{g}}\leq \eta_{\mathrm{n}}$, which translates to the fact that agents performing Newton updates converges at a faster rate. Therefore, as we increase $q$ (more agents perform Newton updates), we obtain faster rate in expectation. This is corroborated in the next section. 

\begin{figure}[t]
  \centering
  \begin{subfigure}[b]{0.49\linewidth}
    \includegraphics[width=\textwidth]{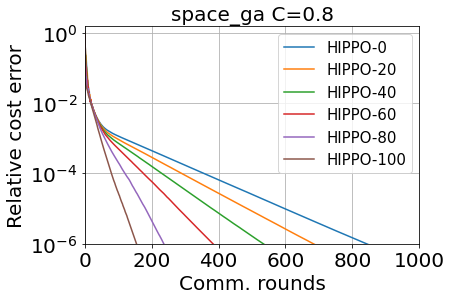}
  \end{subfigure}
  \begin{subfigure}[b]{0.49\linewidth}
    \includegraphics[width=\textwidth]{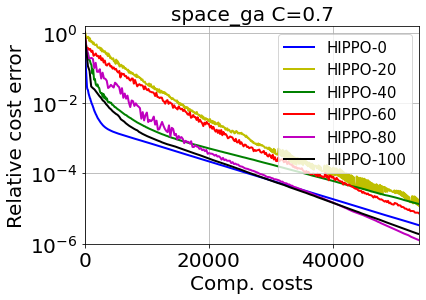}
  \end{subfigure}
  \caption{Performance comparison with different percentage of agents performing Newton updates using comm. rounds (left) and computation costs (right) as the metric.}
  \label{fig1}
\end{figure}
\section{Experiments} 
In this section, we present numerical experiments for HIPPO. We consider the following distributed $\gamma$-weighted LASSO problem:
\begin{gather*}
    \underset{x\in\mathbb{R}^d}{\text{minimize}}\,\,l(x)=\bigg\{ \sum_{i=1}^m \frac{1}{2}\norm{A_ix-b_i}^2+\gamma\norm{x}_1\bigg\}.
\end{gather*}
We generate random graphs with $m=50$ agents by repetitively drawing edges according to a Bernoulli($p$) distribution (with redrawing if necessary to ensure connectivity). We evenly distribute $3,000$ data points from the space\_ga\footnote{https://www.csie.ntu.edu.tw/$\sim$cjlin/libsvm/} dataset to agents, where $d=6$. The averaged relative loss $\frac{\tfrac{1}{m}\sum_{i=1}^m l(x^t_i)-l(x^\star)}{\tfrac{1}{m}\sum_{i=1}^m l(x^0_i)-l(x^\star)}$ is plotted against two metrics, namely communication rounds and computation costs. Note that a beneficial attribute of the proposed method is that the communication costs per agent in one iteration of the algorithm, regardless of the computational choice of the agents. The fraction of active agents is denoted by $C\in(0,1]$, where agents are sampled uniformly at random. We compare the performance of the algorithm with different proportion of agents using Newton updates. For example, HIPPO--20 stands for $20\%$ of agents perform Newton updates and the rest perform gradient updates. We observe that in Fig.\ref{fig1}, as the number of agents performing Newton updates increases, the performance strictly increases. However, when computation costs is used as the metric, the candidate with the optimal performance is not obvious.

\bibliography{ref.bib}
\bibliographystyle{IEEEtran}

\appendix

\textit{Proof of Lemma \ref{lemma1}}: We denote $y^t=[\alpha^t;\beta^t]$ for some $\alpha^t,\beta^t\in\mathbb{R}^{nd}$. By premultiplying (\ref{admm4}) with $B^\top$ and using (\ref{z_update}), it holds that $B^\top y^{t+1}=0$. By the definition of $B$ matrix in (\ref{compact}) and the initialization $\alpha^0=-\beta^0$, we conclude that $\alpha^t=\beta^t$ for $t\geq 0$. Therefore, we can rewrite the dual update (\ref{admm4}) as:
\begin{align}
    \alpha^{t+1} &=\alpha^t+\mu_z (A_s x^{t+1}-z^{t+1}) ,\label{lemma1_proof_1}\\
    -\alpha^{t+1} & = -\alpha^t+\mu_z(A_dx^{t+1}-z^{t+1}).\label{lemma1_proof_2}
\end{align} By taking the sum and the difference of (\ref{lemma1_proof_1}) and (\ref{lemma1_proof_2}), we obtain for $t\geq 0$, 
$
    z^{t+1} = \tfrac{\mz}{2}E_ux^{t+1}, 
    \alpha^{t+1}=\alpha^t+\tfrac{\mz}{2}E_sx^{t+1}.
$
Recalling the definition $\phi^t = E_s^\top \alpha^t$ and by premultiplying the update of $\alpha^{t+1}$ with $E_s^\top$, we obtain the update (\ref{lemma1_dual}). It remains to show that the primal update (\ref{x_primal}) is equivalent to (\ref{lemma1_x}). The primal update (\ref{x_primal}) can be expressed as: 
\begin{gather}
    x^{t+1}=x^t-(H^t)^{-1}\big[\nabla F(x^t)+A^\top y^t+S \lambda^t\nonumber\\
    +\mz A^\top (Ax^t-Bz^t)+\mt S(S^\top x^t-\theta^t)\big]. \label{proof_x_update}
\end{gather}
Since $y^t=[\alpha^t;-\alpha^t]$ and $z^t=\tfrac{1}{2}E_u x^t$, we obtain $A^\top y^t=E_s^\top \alpha^t=\phi^t$ and $\mz A^\top (Ax^t-Bz^t)=\tfrac{\mz}{2}L_ux^t$. After substituting these equalities into (\ref{proof_x_update}), we obtain (\ref{lemma1_x}). \QEDB

\textit{Proof of Lemma 2}: From the primal update (\ref{lemma1_x})and the identities $\phi=E_s^\top \alpha$, $L_s=E_s^\top E_s$, we obtain: 
\begin{gather}
    \nabla F(x^t)+E_s^\top \alpha^t+S\lambda^t+\tfrac{\mu_z}{2}L_s x^t+\mt S (S^\top x^t-\theta^t) \nonumber
    \\
    +H^t(x^{t+1}-x^t)= 0. \label{primal rearrange}
\end{gather}
Using the dual update (\ref{lemma1_dual})-(\ref{lemma1_lambda}), we further obtain: 
\begin{gather}
    E_s^\top \alpha^t+\tfrac{\mu_z}{2}L_sx^t = E_s^\top \alpha^{t+1}-\tfrac{\mu_z}{2}L_s(x^{t+1}-x^t).\label{E_s alpha}\\ 
    S\lambda^{t}+\mt S(S^\top x^t-\theta^t) =S\lambda^{t+1}-\mt S[S^\top(x^{t+1}-x^t)\nonumber\\
    -(\theta^{t+1}-\theta^t)].  \label{lemma2_proof_lambda}
\end{gather}
By substituting (\ref{E_s alpha}) and (\ref{lemma2_proof_lambda}) into (\ref{primal rearrange}), we obtain: 
\begin{gather}
    \nabla F(x^t)+E_s^\top \alpha^{t+1}
    -\tfrac{\mu_z}{2}L_s(x^{t+1}-x^t)+S[\lambda^{t+1}\label{lemma2_proof_33}\\
    -\mt S^\top(x^{t+1}-x^t)
    +\mt (\theta^{t+1}-\theta^t)]+H^t(x^{t+1}-x^t) = 0. \nonumber
\end{gather}
Recall $2D=L_s+L_u$. Adding and subtracting $(\mu_z D+\mt S S^\top+\Delta)(x^{t+1}-x^t)$ from (\ref{lemma2_proof_33}), we obtain:
\begin{gather}
    \nabla F(x^t)+E_s^\top \alpha^{t+1}+\tfrac{\mu_z}{2}L_u(x^{t+1}-x^t)+S\lambda^{t+1}\nonumber\\
    +\mu_\theta S(\theta^{t+1}-\theta^t)
    +\Delta(x^{t+1}-x^t)
    \nonumber\\
    +\left(H^t-\mu_z D-\mt S S^\top-\Delta\right)(x^{t+1}-x^t)=0.\label{lemma2_proof_34}
\end{gather}
From Lemma 1 and KKT conditions, it holds that:  
\begin{gather}
    \tfrac{\mz}{2}L_u(x^{t+1}-x^t)=\mz E_u^\top (z^{t+1}-z^t),\label{lemma2_proof_z}\\
    \nabla F(x^\star)+E_s^\top \alpha^\star+S \lambda^\star = 0.\label{kkta}
\end{gather}
By substituting (\ref{lemma2_proof_z}) into (\ref{lemma2_proof_34}) and subtracting (\ref{kkta}), we obtain the desired. \QEDB

\textit{Proof of Lemma 3}: By the definition of $e^t$ in (\ref{error}): $e^t = \nabla F(x^t)-\nabla F(x^{t+1})+J^t(x^{t+1}-x^t)$, we obtain
\begin{gather*}
    \norm{e^t_i}\leq \norm{\nabla f_i(x_i^t)-\nabla f(x_i^{t+1})}+\norm{J_{ii}^t}\norm{x_i^{t+1}-x_i^t}.
\end{gather*}
For $i\in \mathcal{U}_{\mathrm{Gradient}}^t$, $J_{ii}^t=0$. By assumption \ref{a1}, it holds that:
$
    \norm{e_i^t}\leq M_f \norm{x^{t+1}_i-x_i^t}.
$
For $i\in \mathcal{U}_{\mathrm{Newton}}$, $J^t_{ii}=\nabla^2 f_i(x_i^t)$. Using assumption \ref{a1} again, we obtain:
\begin{gather}
    \norm{e_i^t}\leq 2M_f\norm{x^{t+1}_i-x^t_i}. \label{lemma3_proof_37}
\end{gather}
Moreover, using the fundamental theorem of calculus, we rewrite $\nabla f_i(x_i^{t+1})-\nabla f_i(x_i^t)$ as: 
\begin{gather*}
    \nabla f_i(x_i^{t+1})-\nabla f_i(x+i^t)
    =\\
    \int_0^1 \nabla^2 f_i(sx_i^{t+1}+(1-s)x_i^t)(x_i^{t+1}-x_i^t) ds.
\end{gather*}
After adding and subtracting $\int_0^1 \nabla^2 f_i(x_i^t)(x_i^{t+1}-x_i^t)ds$, we obtain: 
\begin{gather*}
        \nabla f_i(x_i^{t+1})-\nabla f_i(x_i^t)
    =\int_0^1 \nabla^2 f_i(x_i^t)(x_i^{t+1}-x_i^t) ds\\
    +\int_0^1 \left(\nabla^2 f_i(sx_i^{t+1}+(1-s)x_i^t)-\nabla^2 f_i(x_i^t)\right)(x_i^{t+1}-x_i^t) ds.
\end{gather*}
Therefore, we obtain the following chain of inequality: 
\begin{align*}
        &\norm{\nabla f_i(x_i^{t+1})-\nabla f_i(x^t)-\nabla^2 f_i(x_i^t)(x_i^{t+1}-x_i^t)}\\
    &=\norm{\int_0^1 \left(\nabla^2 f_i(sx_i^{t+1}+(1-s)x_i^t)-\nabla^2 f_i(x_i^t)\right)(x_i^{t+1}-x_i^t)ds}\\
    &\leq \int_0^1 \norm{\nabla ^2 f_i(sx_i^{t+1}+(1-s)x_i^t)-\nabla^2 f_i(x_i^t)}\cdot\norm{x_i^{t+1}-x_i^t}ds \\
    &\leq \int_0^1 sL_f\norm{x_i^{t+1}-x_i^t}^2 ds
    =\tfrac{L_f}{2}\norm{x_i^{t+1}-x_i^t}^2. 
\end{align*}
Note that for $i\in \mathcal{U}^t_{\mathrm{Newton}}$, it holds that:
\begin{gather*}
    \norm{e_i^t}=\norm{\nabla f_i(x_i^{t+1})-\nabla f_i(x^t)-\nabla^2 f_i(x_i^t)(x_i^{t+1}-x_i^t)}\\
    \leq \tfrac{L_f}{2}\norm{x_i^{t+1}-x_i^t}. 
\end{gather*}
    
Combining (\ref{lemma3_proof_37}) and the above, we obtain: for $i\in \mathcal{U}^t_{\mathrm{Newton}}$, $\norm{e_i^t}\leq \min\{2M_f,\tfrac{L_f}{2}\norm{x^{t+1}_i-x^t_i}\}\norm{x_i^{t+1}-x^t_i}.$ Using the notation in (\ref{lemma3_eqn}), we obtain the desired. \QEDB

\textit{Proof of Theorem 1}: Part (i): We proceed by establishing that:
\begin{gather}
    \norm{Tv^{t}_\alpha-v^\star_\alpha}^2_{\mathcal{H}}\leq \tfrac{1}{\eta} \norm{v^t_\alpha-v^\star_\alpha}^2_{\mathcal{H}}.\label{theorem1_proof_1}
\end{gather}
We first prove the following useful inequalities:
\begin{align}
       (\lambda^{t+1}-\lambda^t)^\top (\theta^{t+1}-\theta^t) &\geq 0 ,\label{claim1_2}\\
       (\lambda^{t+1}-\lambda^\star)^\top (\theta^{t+1}-\theta^\star) &\geq 0.\label{claim1_1}
\end{align}
From the update (\ref{lemma1_theta}), its optimality condition, and dual update \ref{lemma1_lambda}, we obtain $$
    0\in \partial g(\theta^{t+1})-\mu_\theta(\tfrac{1}{\mu_\theta}\lambda^t+S^\top x^{t+1}-\theta^{t+1})
    =\partial g(\theta^{t+1})-\lambda^{t+1}.
$$ Therefore, we obtain (\ref{claim1_2}) by using assumption \ref{a3} and (\ref{claim1_1}) by using KKT conditions. We proceed to establish a useful upper bound for $\norm{\alpha^{t+1}-\alpha^\star}^2+\norm{\lambda^{t+1}-\lambda^\star}^2$ by first showing that $[\alpha^t;\lambda^t]$ that lies in the column space of $[E_s; S^\top]$. Recalling the dual update (\ref{lemma1_dual})-(\ref{lemma1_lambda}) and the identity $\phi^t=E_s^\top \alpha^t$, it suffices to show that the column space of $[0; \mt I]$ is a subset of the column space of $[\tfrac{\mz}{2}E_s;\mt S^\top]$. Consider fixed $r^x\in\mathbb{R}^d$ and construct a vector $r^y\in\mathbb{R}^{md}$ as $r^y=[r^x;\dots;r^x]$. Then it holds that $$[\tfrac{\mz}{2}E_s;\mt S^\top] r^y=[0;\mt I]r^x,$$ which shows the desired inclusion relation. Moreover, it is easy to verify that there exists a unique pair $[\alpha^\star;\lambda^\star]$ that lies in the column space of $[E_s;S^\top]$. By selecting $\mz=2\mt$, we obtain that $[\alpha^t-\alpha^\star;\lambda^t-\lambda^\star]$ lies in the column space of $[E_s;S^\top]$. Denoting the smallest positive eigenvalue of $[E_s;S^\top] [E_s^\top, S]$ as $\sigma^+_{\mathrm{min}}$, we obtain:
\begin{align}
    &\norm{\alpha^{t+1}-\alpha^\star}^2+\norm{\lambda^{t+1}-\lambda^\star}^2
    \nonumber\\
    &\leq \tfrac{1}{\sigma^+_{\mathrm{min}}}\norm{E_s^\top (\alpha^{t+1}-\alpha^\star)+S(\lambda^{t+1}-\lambda^\star) }^2.\label{lemma_col_space}
\end{align}
Since $F(x)$ is $m_f$-strongly convex and $M_f$-Lipschitz continuous, the following holds:
\begin{align}
    &\tfrac{m_fM_f}{m_f+M_f}\norm{x^{t+1}-x^\star}^2+\tfrac{1}{m_f+M_f}\norm{\nabla F(x^{t+1})-\nabla F(x^\star)}^2 \nonumber\\
    &\leq (x^{t+1}-x^\star)^\top (\nabla F(x^{t+1})-\nabla F(x^\star)) \nonumber \\
    &\leq -(x^{t+1}-x^\star)^\top e^t
    -\epsilon(x^{t+1}-x^\star)^\top(x^{t+1}-x^t)\nonumber\\ &-(x^{t+1}-x^\star)E_s^\top (\alpha^{t+1}-\alpha^\star)\nonumber\\
    &-(x^{t+1}-x^\star)^\top S \left(\lambda^{t+1}-\lambda^\star+\mt(\theta^{t+1}-\theta^t)\right)\nonumber\\
    &-\mz(x^{t+1}-x^\star)^\top E_u^\top (z^{t+1}-z^t),\label{main_1}
\end{align}
where the last inequality follows from Lemma \ref{lemma2}. From the dual update and KKT conditions, we obtain \begin{align*}
        &(x^{t+1}-x^\star)^\top E_s^\top = \tfrac{2}{\mz}(\alpha^{t+1}-\alpha^t)^\top ,\\
    &(x^{t+1}-x^\star)^\top S = (\theta^{t+1}-\theta^\star)^\top +\tfrac{1}{\mt}(\lambda^{t+1}-\lambda^t)^\top,\\
    &(x^{t+1}-x^\star)^\top E_u^\top = 2(z^{t+1}-z^\star)^\top.
\end{align*}
Using these expressions for $(x^{t+1}-x^\star)^\top E_s^\top$,  $(x^{t+1}-x^\star)^\top S$, and $(x^{t+1}-x^\star)^\top E_u^\top$, we rewrite (\ref{main_1}) as:
\begin{align}
    &\tfrac{2m_fM_f}{m_f+M_f}\norm{x^{t+1}-x^\star}^2+\tfrac{2}{m_f+M_f}\norm{\nabla F(x^{t+1})-\nabla F(x^\star)}^2
    \nonumber\\
    &\leq-2(x^{t+1}-x^\star)^\top e^t-2\epsilon (x^{t+1}-x^\star)^\top(x^{t+1}-x^t)\nonumber\\
    &-\tfrac{4}{\mz}(\alpha^{t+1}-\alpha^t)^\top (\alpha^{t+1}-\alpha^\star)
    -\tfrac{2}{\mt}(\lambda^{t+1}-\lambda^t)^\top (\lambda^{t+1}-\lambda^\star)\nonumber\\
    &-2\mt(\theta^{t+1}-\theta^\star)^\top (\theta^{t+1}-\theta^t)
    -4\mz(z^{t+1}-z^\star)^\top(z^{t+1}-z^t)\nonumber\\
    &=\norm{\va^t-\vas}^2_{\mathcal{H}}-\norm{\va^{t+1}-\vas}^2_\mathcal{H}-\norm{\va^{t+1}-\va^t}^2_{\mathcal{H}}\nonumber\\
    &-2(x^{t+1}-x^\star)^\top e^t, \label{theorem1_proof_45}
\end{align}
where have used (\ref{claim1_1}), (\ref{claim1_2}), and the identity: $-2(a-b)^\top  (a-c)=\norm{b-c}^2-\norm{a-b}^2-\norm{a-c}^2$. To establish (\ref{theorem1_proof_1}), we need to show for some $\eta>0$,
\begin{align}
    \eta\norm{v^{t+1}_\alpha-v^\star_\alpha}^2_{\mathcal{H}}\leq \norm{v^t_\alpha-v^\star_\alpha}^2_{\mathcal{H}}-\norm{v^{t+1}-v^\star_\alpha}^2_{\mathcal{H}}. \label{theorem1_proof_linear}
\end{align}
Using (\ref{theorem1_proof_45}), it suffices to show:
\begin{gather}
    \eta \norm{v^{t+1}_\alpha-v^\star_\alpha}^2_{\mathcal{H}}\leq \tfrac{2m_fM_f}{m_f+M_f}\norm{x^{t+1}-x^\star}^2\label{theorem1_proof_suffice}\\
    +\tfrac{2}{m_f+M_f}\norm{\nabla F(x^{t+1})-\nabla F(x^\star)}^2+2(x^{t+1}-x^\star)^\top e^t.\nonumber
\end{gather}
We expand $\eta \norm{v^{t+1}_\alpha-v^\star_\alpha}^2_{\mathcal{H}}$ as follows:
\begin{gather*}
     \eta\norm{v^{t+1}_\alpha-v^\star_\alpha}^2_{\mathcal{H}}= \eta\bigg(\epsilon\norm{x^{t+1}-x^\star}^2
     +2\mz\norm{z^{t+1}-z^\star}^2\\
    +\tfrac{2}{\mz}\norm{\alpha^{t+1}-\alpha^\star}^2 
    +\mt\norm{\theta^{t+1}-\theta^\star}^2
    +\tfrac{1}{\mt}\norm{\lambda^{t+1}-\lambda^\star}^2\bigg).
\end{gather*}
We proceed to establish upper bounds for each component. From Lemma 2, the following holds:
\begin{gather*}
     E_s^\top(\alpha^{t+1}-\alpha^\star)+S (\lambda^{t+1}-\lambda^\star)
    =-\big\{\nabla F(x^{t+1})-\nabla F(x^\star)\nonumber\\
    +\epsilon(x^{t+1}-x^t)
    +\mz E_u^\top(z^{t+1}-z^t)+\mt S(\theta^{t+1}-\theta^t)+e^t\big\}.
\end{gather*}
Using (\ref{lemma_col_space}) and selecting $\mz=2\mt$, we further obtain:
\begin{align}
    &\tfrac{2}{\mz}\norm{\alpha^{t+1}-\alpha^\star}^2+\tfrac{1}{\mt} \norm{\lambda^{t+1}-\lambda^\star}^2 
    \leq \tfrac{5}{\mt\sigma^+_{\mathrm{min}}}(\norm{e^t}^2\nonumber\\
    &+\norm{\nabla F(x^{t+1})-\nabla F(x^\star)}^2+\epsilon^2\norm{x^{t+1}-x^t}^2\nonumber\nonumber\\
    &+\mt^2\norm{\theta^{t+1}-\theta^t}^2+\sigma^{L_u}_{\mathrm{max}}\mz^2\norm{z^{t+1}-z^t}^2),\label{theorem_linear_part1}
\end{align}
where $\sigma^{L_u}_{\mathrm{max}}$ is the maximum eigenvalue of $L_u=E_u^\top E_u$. Moreover, since $z^{t+1}-z^\star= \tfrac{1}{2}E_u(x^{t+1}-x^\star)$, it holds that
$
    2\mz\norm{z^{t+1}-z^\star}^2\leq \tfrac{\mz\sigma^{L_u}_{\mathrm{max}}}{2}\norm{x^{t+1}-x^\star}^2.
$
Similarly, from the dual update (\ref{lemma1_theta}), the following holds: 
$
    \mt\norm{\theta^{t+1}-\theta^\star}^2\leq 2\mt \norm{x^{t+1}-x^\star}^2+\tfrac{2}{\mt} \norm{\lambda^{t+1}-\lambda^t}^2.
$
Using these bounds and considering (\ref{theorem1_proof_suffice}), we need to show the following holds for $\eta>0$,
\begin{gather}
    \eta \norm{\va^{t+1}-\vas}^2_{\mathcal{H}}\leq \eta\bigg\{\tfrac{5}{\mt\sigma^+_{\mathrm{min}}}\Big[\norm{\nabla F(x^{t+1})-\nabla F(x^\star)}^2\nonumber\\
    +\epsilon^2\norm{x^{t+1}-x^t}^2
    +\mt^2\norm{\theta^{t+1}-\theta^t}^2
    +\norm{e^t}^2\nonumber\\
    +\sigma^{L_u}_{\mathrm{max}}\mz^2\norm{z^{t+1}-z^t}^2\Big]
    +\tfrac{2}{\mt}\norm{\lambda^{t+1}-\lambda^t}^2\nonumber\\
    +(\epsilon+2\mt+\tfrac{\mz\sigma^{L_u}_{\mathrm{max}}}{2})\norm{x^{t+1}-x^\star}^2\bigg\}\leq\nonumber\\
     \tfrac{2m_fM_f}{m_f+M_f}\norm{x^{t+1}-x^\star}^2+\tfrac{2}{m_f+M_f}\norm{\nabla F(x^{t+1})-\nabla F(x^\star)}^2\nonumber\\
    +\norm{\va^{t+1}-\va^t}^2_{\mathcal{H}}+2(x^{t+1}-x^\star)^\top e^t.
\end{gather}
Note that $-\zeta \norm{e^t}^2-\tfrac{1}{\zeta}\norm{x^{t+1}-x^\star}^2\leq 2(x^{t+1}-x^\star)^\top e^t$ holds for any $\zeta>0$, it is sufficient to show the following to establish (\ref{theorem1_proof_linear}):
\begin{align*}
     &\zeta\norm{x^{t+1}-x^t}^2_{(\Pi^t)^2}+\eta\bigg\{\tfrac{5}{\mt\sigma^+_{\mathrm{min}}}\Big[\norm{\nabla F(x^{t+1})-\nabla F(x^\star)}^2\\
      &+\norm{x^{t+1}-x^t}^2_{\epsilon^2 I+(\Pi^t)^2}
    +\mt^2\norm{\theta^{t+1}-\theta^t}^2\\
    &+\sigma^{L_u}_{\mathrm{max}}\mz^2\norm{z^{t+1}-z^t}^2\Big]+\tfrac{2}{\mt}\norm{\lambda^{t+1}-\lambda^t}^2\\
    &+(\epsilon+2\mt+\tfrac{\mz\sigma^{L_u}_{\mathrm{max}}}{2})\norm{x^{t+1}-x^\star}^2\bigg\}
    \leq\\
    &\left(\tfrac{2m_fM_f}{m_f+M_f}-\tfrac{1}{\zeta}\right)\norm{x^{t+1}-x^\star}^2
    +\epsilon\norm{x^{t+1}-x^t}^2\\
    &+2\mz\norm{z^{t+1}-z^t}^2
    +\tfrac{2}{\mz}\norm{\alpha^{t+1}-\alpha^t}^2+\mt\norm{\theta^{t+1}-\theta^t}^2\\
    &+\tfrac{1}{\mt}\norm{\lambda^{t+1}-\lambda^t}^2
    +\tfrac{2}{m_f+M_f}\norm{\nabla F(x^{t+1})-\nabla F(x^\star)}^2. \end{align*}
By selecting $\eta$ as in (\ref{eta}), we obtain the desired. Moreover observe that:
\begin{gather}
    \eta \leq \tfrac{\mt \sigma^+_{\mathrm{min}}}{5}\cdot\tfrac{\norm{x^{t+1}-x^t}^2_{\epsilon I-\zeta(\Pi^t)^2}}{\norm{x^{t+1}-x^t}^2_{\epsilon^2 I+(\Pi^t)^2}} \label{hybrid_eta} .
\end{gather}
We have now established (\ref{theorem1_proof_1}). We proceed to show the asynchronous updates converges linearly in expectation. Note that: 
\begin{align*}
    &\norm{\va^{t+1}-\vas}^2_{\mathcal{H}\Omega_\alpha^{-1}}
    = \norm{\va^{t}+\Omega_\alpha^{t+1}(T\va^t-\va^t)-\vas}^2_{\mathcal{H}\Omega_\alpha^{-1}} \\
    &= \norm{\va^t-\vas}_{\mathcal{H}\Omega_\alpha^{-1}}^2 +2(\va^t-\vas)^\top\mathcal{H}\Omega_\alpha^{-1}\Omega_\alpha^{t+1}(T\va^t-\va^t)\\
     &+ (T\va^t-\va^t)^\top\Omega_\alpha^{t+1}\mathcal{H}\Omega_\alpha^{-1}\Omega_\alpha^{t+1}(T\va^t-\va^t).
\end{align*}
Since $\Omega_\alpha^{t+1},\Omega^{-1}_\alpha$, and $\mathcal{H}$ are all diagonal matrices, they commute with each other. Moreover, since each sub-block of $\Omega_\alpha^{t+1}$ is $I_d$ or 0, it holds that $\Omega_\alpha^{t+1}\Omega_\alpha^{t+1}=\Omega^{t+1}_\alpha$. Therefore, 
\begin{align*}
    &\mathbb{E}^t\left[\norm{\va^{t+1}-\vas}^2_{\mathcal{H}\Omega_\alpha^{-1}} \right] 
     = \norm{\va^{t}-\vas}_{\mathcal{H}\Omega_\alpha^{-1}}^2
+\norm{T\va^t-\va^t}^2_{\mathcal{H}}\\
&+2(\va^t-\vas)^\top\mathcal{H}(T\va^t-\va^t) \\
&\leq \norm{\va^t-\vas}_{\mathcal{H}\Omega_\alpha^{-1}}^2-\tfrac{\eta}{1+\eta}\norm{\va^t-\vas}^2_{\mathcal{H}} \\
&\leq  \left(1-\tfrac{p^{\mathrm{min}}\eta }{1+\eta} \right)\norm{\va^t-\vas}^2_{\mathcal{H}\Omega_\alpha^{-1}},
\end{align*}
which is the desired. 
\QEDB

\textit{Proof of Corollary 1}: We first distribute each $\alpha_k,k\in[n]$, to each edge and label agents and edges with an arbitrary order. For each edge $\mathcal{E}_k$, we write $\mathcal{E}_k=(i,j)$ with the convention $i<j$. Then for each agent $i$, we can divide the incident edges to two groups:  $\mathcal{P}_i=\{k:\mathcal{E}_k=(i,j),j\in\mathcal{N}_i\}$ and $\mathcal{S}_i=\{k:\mathcal{E}_k=(j,i),j\in\mathcal{N}_i\}$. Consider the activation scheme using $\Omega^{t+1}$. Recall $\alpha^{t+1}_k=\alpha^t_k+\tfrac{\mz}{2}(x_i^{t+1}-x_j^{t+1})$. Then the dual updates are described by:
\begin{align*}
    &\phi^{t+1}_i = \phi^t_i +\tfrac{\mz}{2}X^{t+1}_{ii}\sum_{j\in\mathcal{N}_i}(x_i^{t+1}-x_j^{t+1})\\
    &= \phi_i^t +X^{t+1}_{ii} \left\{\sum_{k\in \mathcal{P}_i}(\alpha^{t+1}_k-\alpha^{t}_k)+\sum_{k\in \mathcal{S}_i}(\alpha^{t}_k-\alpha^{t+1}_k)\right\}.
\end{align*}
Therefore, if $X^{t+1}_{ii}=I_d$, then $Y^{t+1}_{kk}=I_d$ for $k\in P_i\cup S_i$ for the corresponding $\Omega^{t+1}_\alpha$, i.e., all incident edges are active. It can be verified that we can map $X^{t+1}$ to $Y^{t+1}$ as:
$$
    Y^{t+1} = \textbf{Blkdiag}\left(\ceil{\frac{E_u X^{t+1} (\mathbf{1}\otimes I_d)}{2}} \right),
$$
where $\lceil\cdot \rceil$ is the entry-wise ceiling operation and $\mathbf{1}\in\mathbb{R}^{m}$ is the all one vector. To show $\mathbb{E}^t[\Omega^{t+1}_\alpha]\succ0$, we only need to show $\mathbb{E}^t[Y^{t+1}]\succ 0$, which amounts to showing that $\mathbb{E}^t\left[\ceil{\frac{E_u X^{t+1} (\mathbf{1}\otimes I_d)}{2}}_{k}\right]\in \mathbb{R}^{d\times d},k\in[n],$ is positive definite. We can explicitly express each block as:
$$
    \ceil{\frac{E_u X^{t+1} (\mathbf{1}\otimes I_d)}{2}}_{k} = \ceil{\frac{X^{t+1}_{ii}+X^{t+1}_{jj}}{2}},
$$
where $(i,j)\in \mathcal{E}_k$. Therefore, 
$$
    \mathbb{E}^t\left[\ceil{\frac{E_u X^{t+1} (\mathbf{1}\otimes I_d)}{2}}_{k}\right] = \mathbb{E}^t \left[\ceil{\frac{X^{t+1}_{ii}+X^{t+1}_{jj}}{2}}\right]\succ 0,
$$
which shows that $\mathbb{E}^t[Y^{t+1}]\succ 0$.\QEDB

%%%%%%%%%%%%%%%%%%%%%%%%%%%%%%%%%%%%%%%%%%%%%%%%%%%%%%%%%%%%%%%%%%%%%%%%%%%%%%%%

\end{document}